\documentclass{amsart}
\usepackage{amssymb}
\usepackage{amsfonts}
\usepackage{amsmath}
\usepackage[legalpaper,bookmarks=true,colorlinks=true,linkcolor=blue,citecolor=blue]{hyperref}
\usepackage{graphicx}%
\usepackage{fancyhdr}
\usepackage{color}
\usepackage[mathlines]{lineno}
\usepackage{lscape}
\usepackage{epsfig}
\usepackage{natbib}
\usepackage{geometry}
\usepackage{tgbonum}
\usepackage[]{listings}

\fontfamily{qcr}\selectfont

\newtheorem{theorem}{Theorem}
\theoremstyle{plain}

\newtheorem{corollary}{Corollary}

\newtheorem{definition}{Definition}

\newtheorem{proposition}{Proposition}

\numberwithin{equation}{section}

\newcommand{\Bin}{\bigskip \noindent}
\newcommand{\Bi}{\bigskip}
\newcommand{\Ni}{\noindent}

\begin{document}
\Large
\title[Asymptotics of summands of associated random variables]{$\mathbf{G}$-Central limit theorems and $\mathbf{G}$-invariance principles for associated random variables}
\author{Aladji Babacar Niang} 
\author{Akym Adekpedjou}
\author{Harouna Sangar\'e}
\author{Gane Samb Lo}

\begin{abstract} The investigation asymptotic limits on associated data mainly focused on limit theorems of summands of associated data and on the related invariance principles. In a series of papers, we are going to set the general frame of the theory by considering an arbitrary infinitely decomposable (divisible) limit law for summands and study the associated functional laws converging to L\'evy processes. The asymptotic frame of Newman (1980) is still used as a main tool. Detailed results are given when $G$ is a Gaussian law (as confirmation of known results) and when $G$ is a Poisson law. In the later case, classical results for independent and identically distributed data are extended to stationary and non-stationary associated data.\\   

\noindent $^{\dag}$ Aladji Babacar Niang\\
LERSTAD, Gaston Berger University, Saint-Louis, S\'en\'egal.\\
Email: niang.aladji-babacar@ugb.edu.sn, aladjibacar93@gmail.com\\

\Ni $^{\dag \dag}$ Akim Adekpedjou, PhD\\
Professor of Statistics\\
Department of Mathematics \& Statistics\\
Missouri University of Science \& Technology\\
Office: +1(573)-341-4649\\
E-mail: akima@mst.edu\\
URL: web.mst.edu/~akima\\
ROlla, MO 65409\\

\noindent $^{\dag \dag \dag}$ Harouna Sangar\'e\\
Main Affiliation: \textit{DER} \textit{MI}, \textit{FST}, Universit\'e des Sciences, des Techniques et des Technologies de Bamako (\textit{USTT-B}), Mali.\\
Affiliation : \textit{LERSTAD}, Universit\'e Gaston Berger (\textit{UGB}), Saint-Louis, S\'en\'egal.\\
Email : harounasangare@fst-usttb-edu.ml, harouna.sangare@mesrs.ml\\
sangare.harouna@ugb.edu.sn, harounasangareusttb@gmail.com\\

\noindent $^{\dag \dag \dag \dag}$ Gane Samb Lo.\\
\textit{LERSTAD}, Gaston Berger University, Saint-Louis, S\'en\'egal (main affiliation).\newline
\textit{LSTA}, Pierre and Marie Curie University, Paris VI, France.\newline
\textit{AUST} - African University of Sciences and Technology, Abuja, Nigeria\\
gane-samb.lo@edu.ugb.sn, gslo@aust.edu.ng, ganesamblo@ganesamblo.net\\
Permanent address : $1178$ Evanston Dr \textit{NW} $T3P 0J9$, Calgary, Alberta, Canada.\\

\textbf{\textit{AMS} $2010$ Mathematics Subject Classification:} $60GXX$, $62GXX$ \\
\end{abstract}
\maketitle

\section{Introduction} \label{sec1}

\Ni Recently, the investigation of fundamental aspects of probability theory on associated data has been very active. Since the pioneering work of \cite{newman}, the type of dependence called association is highly present in the literature. Since independent sequences are associated, the research community tries to extend of results of the independence frame to the association frame. Because this concept is not popular yet, we are going to give a brief account of it below, before we describe our aims and motivation with more precision.\\

\subsection{A brief account on associated data} $ $\\

\noindent The definition is given by \cite{esary} as follows.

\begin{definition} \label{def-association}
A finite sequence of \textit{rv}'s $(X_{1},\cdots,X_{n})$ are associated when for any couple of real valued and coordinate-wise non-decreasing functions $h$ and $g$ defined on $\mathbb{R}^{n}$, we have 

\begin{equation}
\mathbb{C}ov(h(X_{1},\cdots,X_{n}),\ \ g(X_{1},\cdots,X_{n}))\geq 0.
\end{equation}%

\Bin An infinite sequence of \textit{rv}'s are associated whenever all its finite
sub-sequences are associated.
\end{definition}

\Bin There is a few number of interesting properties to be found in \cite{rao} :

\bigskip

\noindent \textbf{(P1)} A sequence of independent \textit{rv}'s is associated. \textbf{(P2)} Partial sums of associated \textit{rv}'s are associated. \textbf{(P3)} Order statistics of independent \textit{rv}'s are associated. \textbf{(P4)} Non-decreasing functions and non-increasing functions of associated variables are associated.  \textbf{(P5)} Suppose that $Z_{1},Z_{2},\cdots,Z_{n}$ is
associated and  $(a_{i})_{1\leq i\leq n}$ is sequence of non-negative  numbers and $(b_{i})_{1\leq i\leq n}$ real numbers, then the \textit{rv}'s $a_{i}(Z_{i}-b_{i})$ are associated.\newline

\Ni The concept of association also arises in Physics, and is quoted under the name of \textit{FKG} property (\cite{fortuin}), in percolation theory and even in Finance (see \cite{jiazhu}).\\

\Ni However, the concept took its place in pure probability theory only with papers by \cite{newmanwright,newmanwright2} which addressed the central limit theorem, the strong law of large numbers and the functional limiting laws. Since then, a great number of contributions appeared in the literature: [\cite{newman}, \cite{sang15}, etc., for example for the strong law of large numbers], [\cite{yu93}, \cite{GH-louhichi}, \cite{sang20}, etc, for the empirical processes;  \cite{cox}, \cite{GH-paulo}, \cite{loSanMbarka18}, \cite{sanglo28}, \cite{akim21} for the central limit theorem, \cite{GH-paulo} for the law of the iterated logarithm, \cite{burton}, \cite{dabro}, \cite{newmanwright2}, \cite{traore16} for the invariance principles, etc.]. There are three important manuscripts aiming at reviewing the current research [\cite{GH-bulinski2007}, \cite{GH-paulo}, \cite{rao}]. \cite{sang15} provides a short but interesting review. Important results in mathematical statistics using associated data also exist and some of them are given in \cite{rao} in particular.\\

\Ni From the review above, wee see that the central limit theorem and the invariance principle topics are central in the asymptotic theory on associated data. However, it appeared that so far the investigations focused in the case where the weak limits are Gaussian random variables [for the \textit{CLT}] or re-scaled Brownian motions [for the invariance principles]. This motivates us to have a general theory for \textit{CLT}'s with arbitrary limiting law $G$ for an infinitely decomposable random variable corresponding to invariance principles to L\'evy processes. This extension is really interesting since the scope of the study passes from one limiting law to an infinite number of limiting laws. From the general results we obtain, we particularize Poisson limits (for \textit{CLT}'s) and compound Poisson processes (for invariance principles).\\

\Ni We organize the rest of the paper as follows. We will present the main notation in Subsection \ref{subsection11} where we recall the most important results on limits of summands for independent random variables with finite variances. Therein, we also describe the general method of using the Newman approximation method called the \textit{block methods} and state intermediate results on the convergence of rates induced by Newman method. In Section \ref{sec2}, we use the general results of Subsection \ref{subsection15} of Section \ref{sec1} to have a general theory of the \textit{CLT} to associated data for both sequences and arrays, for both stationary and non-stationary. In Section \ref{sec3}, we specialized the results in Subsection \ref{subsection15} for stationary associated sequences. In section \ref{sec4}, we consider two iconic examples of convergence to Poisson law, namely the summands of Bernoulli \textit{rv}.'s and corrected geometric \textit{rv}.'s in the by-row \textit{iid} case, and extend the related results to non-stationary by-row independent arrays and next to the association frame for both stationary and non-stationary.\\

\Ni The obtained \textit{CLT}'s should be extended to invariance principles. However in fear of presenting a too lengthy paper, we decided to devote a coming paper to those invariance principles.\\

\subsection{Notation and $G$-\textit{CLT} for summands of random variables random variables} \label{subsection11}$ $\\

\Ni Let us consider the array
  
$$
X\equiv \biggr\{ \{X_{k,n}, \ 1 \leq k \leq k_n=k(n)\}, \ n\geq 1\biggr\},
$$

\Bin of square integrable random variables defined on the same probability space $(\Omega, \mathcal{A}, \mathbb{P})$. We denote $F_{k,n}$ as the cumulative distribution function (\textit{cdf}) of $X_{k,n}$. We also denote by $a_{k,n} = \mathbb{E}(X_{k,n})$ and $\sigma_{k,n}^2=\mathbb{V}ar(X_{k,n})$, $1\leq k \leq k(n)$, if these expectations or variances exist. We also suppose that 

\begin{equation} \label{C0}
k(n) \rightarrow +\infty \  as \ n\rightarrow +\infty. 
\end{equation}

\Bin The central limit theorem problem consists in finding, whenever possible, the weak limit law of the by-row sums of the array $X$,  \textit{i.e.} the summands:

$$
S_n[X]=\sum_{k=1}^{k(n)} X_{k,n}, \ n\geq 1.
$$ 

\Bin Historically, the \textit{CLT} was discovered with the convergence of a Binomial law (which has the same law as the sum of \textit{iid} Bernoulli random variables) to the standard Gaussian law (due to Laplace, De Moivre, etc., around $1731$, see \cite{loeve} for a review). For a long period of time, the Gaussian limit was automatically meant in the \textit{CLT} problem. Many authors, among them L\'evy, Gnedenko, Kolmogorov, etc., characterized the class of possible limit laws under the \textit{Uniform Asymptotic Negligibility (\textit{UAN})} condition, exactly as the class of infinitely decomposable distributions (see below). The longtime association of \textit{CLT}'s to Gaussian limit explains that some authors reserve the vocable \textit{CLT} for the Gaussian limit and for other possible limits, they use the vocable \textbf{Non-Central limit theorem}. Here we use $\textbf{$G$-\textit{CLT}}$ to cover all possible limit laws $G$ beyond the Gaussian law.\\

\Ni Here we suppose that the $X_{k,n}$'s are integrable with finite variances. For an array $X$, we define some important hypotheses used in the formulation of the \textit{CLT} problem.\\

\Ni \textbf{(1) The \textit{UAN} condition}: for any $\varepsilon>0$,

\begin{equation}
U(n,\varepsilon,X)= \sup_{1\leq k \leq k(n)} \mathbb{P}(|X_{k,n} - a_{k,n}|\geq \varepsilon)\rightarrow 0.
\end{equation}

\Bin \textbf{(2) The Bounded Variance Hypothesis (\textit{BVH})}: there exists a constant $c>0$,

$$
\sup_{n\geq 1} MV(n,X) \leq c,
$$

\Bin where

$$
MV(n,X)=\mathbb{V}ar(S_n[X]), \ n\geq 1.
$$

\Bin \textbf{(3) The Variance Convergence Hypothesis (\textit{VCH})}:

$$
MV(n,X) \rightarrow c \in ]0,+\infty[.
$$

\Bin According to the state of the art in \textit{CLT}'s theory for centered, square integrable  and independent array of random variables, the summands weakly converges to a \textit{rv} associated to the \textit{cdf} $G$ and to the \textit{cha.f} $\psi_{G}$ under the \textit{UAN} condition and the \textit{BVH} if and only if the sequence of distribution functions (\textit{df})

$$
K_n(x)=\sum_{k=1}^{k(n)} \int_{-\infty}^{x} y^2 dF_{k,n}(y), \ x\in \mathbb{R}, \ n\geq 1,
$$

\Bin pre-weakly converges to a \textit{df} $K$, denoted $K_n \rightsquigarrow_{pre} K$, that is for any continuity point $x$ of $K$ denoted as $[x \in C(K)]$,

$$
K_n(x) \rightarrow K(x),
$$

\Bin and the \textit{cha.f} $\psi_{G}(\circ)$ of $G$ is given by $\exp(\psi[K](\circ))$ with

\begin{equation}
\forall u \in \mathbb{R}, \ \psi[K](u)=\int \frac{e^{iux}-1-iux}{x^2} \ dK(x).
\end{equation}

\Bin If we have the \textit{VCH}, the convergence criterion is replaced by the weak convergence $K_n \rightsquigarrow K$. Moreover, the limit law $G$ is necessarily an infinitely decomposable law.\\

\Ni In the non centered case, with the same hypotheses above on the random variables of the array, the summands weakly converges to a probability law  associated to the \textit{cdf} $G^{\ast}$ and to the \textit{cha.f} $\psi_{G^{\ast}}$ under the \textit{UAN} condition and the \textit{BVH} if and only if 

$$
\sum_{k=1}^{k(n)} a_{k,n} \rightarrow a, \ a\in\mathbb{R}
$$

\Bin and the sequence of distribution functions (\textit{df})

$$
K_n^{\ast}(x)=\sum_{k=1}^{k(n)} \int_{-\infty}^{x} y^2 dF_{k,n}(y+a_{k,n}), \ x\in \mathbb{R}, \ n\geq 1,
$$

\Bin pre-weakly converges to a \textit{df} $K^{\ast}$ and the \textit{cha.f} $\psi_{G^{\ast}}(\circ)$ of $G^{\ast}$ is given by $\exp(\psi[K^{\ast}](\circ))$ with

\begin{equation*}
\forall u \in \mathbb{R}, \ \psi[K^{\ast}](u)=\int \frac{e^{iux}-1-iux}{x^2} \ dK^{\ast}(x).
\end{equation*}

\Bin If we have the \textit{VCH}, the convergence criterion is replaced by the weak convergence $K_n^{\ast} \rightsquigarrow K^{\ast}$. Moreover, the limit law $G^{\ast}$ is of the form $G^{\ast} = G+a$, with $G$ is necessarily a centered and infinitely decomposable law.\\

\subsection{The data-regrouping method} $ $\\

\Ni In many situations with dependent data, there are efforts to transform the summands $S_n[X]$ into blocks with increasing lengths such that the dependence between the summands of blocks decreases as $n$ grows to infinity. We are going to apply that idea here.\\

\Ni We consider sequences of non-negative integer numbers $(m(n),\ \ell(n), \ r(n))$ such that:

$$
\forall n\geq 1, \  k(n)=m(n) \ell(n) + r(n)
$$

\Bin and

$$
1\leq m(n) \rightarrow +\infty, \ \ 0\leq r(n)<\ell(n) \ \ and \ \ \ell(n)/m(n) \rightarrow 0.
$$  

\Bin Accordingly, we define

$$
Y_{j,n}=\sum_{h=1}^{\ell(n)} X_{(j-1)\ell(n)+h,n}, \ 1\leq j\leq m(n); \ \ Y_{n}^{\ast}=\sum_{h=1}^{r(n)} X_{m(n)\ell(n)+h,n} 
$$

\Bin and

$$
S_{m(n)\ell(n)}=\sum_{j=1}^{m(n)\ell(n)} X_{j,n} = \sum_{j=1}^{m(n)} Y_{j,n}.
$$

\Bin So, we have

$$
S_{n}=\sum_{j=1}^{m(n)\ell(n)} X_{j,n} + Y_{n}^{\ast}= S_{m(n)\ell(n)} + Y_{n}^{\ast}, \ n\geq 1.
$$

\Bin We denote the new array as

$$
Y\equiv \biggr\{ \{Y_{j,n}, \ 1 \leq j \leq m(n)\}, \ n\geq 1\biggr\}.
$$

\Bin We can denote 

$$
S_{m(n)}[Y] = \sum_{j=1}^{m(n)} Y_{j,n}, \ \ s_{j,n}^{2} = \mathbb{V}ar(Y_{j,n}) \ \ and \ \ s_{\ast,n}^{2}=\mathbb{V}ar(Y_{n}^{\ast}).
$$

\Bi

\subsection{A general method for handling the blocks} $ $\\

\Bin We need two steps:\\

\Ni \textbf{Step 1}. We try to transform the study of the summands $S_n[X]$ by those of $S_{m(n)}[Y]$ by getting rid of $Y_{n}^{\ast}$. It is easy to check that for any $u\in \mathbb{R}$,

\begin{equation} \label{C1}
\left|\psi_{S_n[X]}(u)-\psi_{S_{m(n)}[Y]}(u) \right| \leq \mathbb{V}ar(Y_{n}^{\ast})^{1/2} \rightarrow 0. 
\end{equation}
 
\Bin \textbf{Step 2}. If \textbf{step 1} is validated, we try to replace the summands $S_{m(n)}[Y]$ by the summands $S_{m(n)}[T]$ where the array

$$
T\equiv \biggr\{ \{T_{j,n}, \ 1 \leq j \leq m(n)\}, \ n\geq 1\biggr\},
$$

\Bin has independent elements in each row and

$$
\forall n\geq 1, \ \forall 1\leq j \leq m(n), \ T_{j,n}=_{d} Y_{j,n}.
$$

\Bin We need to have for any $u\in \mathbb{R}$,

\begin{equation} \label{C2}
\left| \psi_{S_{m(n)}[Y]}(u)-\psi_{S_{m(n)}[T]}(u) \right| \rightarrow 0, 
\end{equation}

\Bin that is

$$
\left| \psi_{S_{m(n)}[Y]}(u)-\prod_{j=1}^{m(n)}  \psi_{Y_{j,n}}(u) \right| \rightarrow 0. 
$$

\Bin Now let us apply the methodology to associated data. \\

\subsection{How works the block method for associated data?} \label{subsection15} $ $\\
\label{description}
\Ni We begin by requiring (\ref{C1}) as an hypothesis and we recall that it holds whenever each row, \textit{i.e.} the $(X_{k,n})_{1\leq k \leq k(n)}$, is a stationary sequence.\\

\Ni Let us address (\ref{C2}). Let us suppose that the elements of the $k(n)$-th row $(X_{k,n})_{1\leq k \leq k(n)}$ are associated, so are the $(Y_{j,n})_{1\leq j \leq m(n)}$. By a Newman's Lemma (See \cite{GH-newmanwright,GH-newmanwright2}), we have for any $u\in \mathbb{R},$

\begin{eqnarray}
\left| \psi_{S_{m(n)}[Y]}(u)-\psi_{S_{m(n)}[T]}(u) \right| \leq R_{m(n),\ell(n)}(u),
\end{eqnarray}

\Bin with 

$$
R_{m(n),\ell(n)}(u)=\frac{u^2}{2}\sum_{1\leq r \neq s\leq m(n)} \mathbb{C}ov(Y_{r,n}, \ Y_{s,n}).
$$

\Bin So Hypothesis (\ref{C2}) is achieved through

\begin{equation} \label{cond_suff_c2}
R_{m(n),\ell(n)}(u)\rightarrow 0. 
\end{equation}

\Bin We get the following general theorem. \\

\begin{theorem} \label{main_theo1} (General Theorem for centered random variables) Let 

$$
X\equiv \biggr\{ \{X_{k,n}, \ 1 \leq k \leq k_n=k(n)\}, \ n\geq 1\biggr\},
$$

\Bin be an array of centered, square integrable and associated random variables defined on the same probability space $(\Omega, \mathcal{A}, \mathbb{P})$. \\

\Ni We denote $F_{k,n}$ as the cumulative distribution function (\textit{cdf}) of $X_{k,n}$ and let

$$
S_n[X]=\sum_{k=1}^{k(n)} X_{k,n} \ \ and \ \ s_n^2[X]=\mathbb{V}ar(S_n[X]), \ n\geq 1.
$$

\Bin Let $(m(n),\ \ell(n), \ r(n))_{n\geq 1}$ be sequences of non-negative integer numbers such that $k(n)=m(n) \ell(n) + r(n)$ for any $n\geq 1$ and 

$$
1\leq m(n) \rightarrow +\infty, \ \ 0\leq r(n)<\ell(n) \ \ and \ \ \ell(n)/m(n) \rightarrow 0.
$$  

\Bin We define 

$$
Y_{j,n}=\sum_{h=1}^{\ell(n)} X_{(j-1)\ell(n)+h,n}, \ 1\leq j\leq m(n); \ \ Y_{n}^{\ast}=\sum_{h=1}^{r(n)} X_{m(n)\ell(n)+h,n} 
$$

\Bin and an array

$$
T\equiv \biggr\{ \{T_{j,n}, \ 1 \leq j \leq m(n)\}, \ n\geq 1\biggr\},
$$

\Bin which has independent elements in each row and

$$
\forall n\geq 1, \ \forall 1\leq j \leq m(n), \ T_{j,n}=_{d} Y_{j,n}.
$$

\Bin We set 

$$
S_{m(n)}[T]=\sum_{j=1}^{m(n)} T_{j,n}
$$

\Bin and denote $G_{j,n}$ the \textit{cdf} of $T_{j,n}$. Then, 
given (\ref{C0}), (\ref{C1}), (\ref{C2}) and assuming that \textit{UAN[T]}, \textit{BVH[T]} [or \textit{CVH[T]}] hold, we have the following characterizations.\\

\Ni (1) $S_n[X]$ converges to a random variable with \textit{cdf} $G$ if and only if $G$ is \textit{idecomp}. \\

\Ni (2) $S_n[X]$  converges to a random variable with \textit{df} $G$ if and only if $\psi_{G}(\circ)=\exp(\psi[K](\circ))$, where

\begin{equation}
\forall u \in \mathbb{R}, \ \psi[K](u)=\int \frac{e^{iux}-1-iux}{x^2} \ dK(x),
\end{equation}

\Bin with $K$ is \textit{df} on $\mathbb{R}$ such that

\begin{equation}
K_n(\circ)= \sum_{j=1}^{m(n)} \int_{-\infty}^{\circ} y^2 \ dG_{j,n}(y) \rightsquigarrow_{pre} K,
\end{equation} 

\Bin under \textit{BVH[T]} and 

\begin{equation}
K_n(\circ)= \sum_{j=1}^{m(n)} \int_{-\infty}^{\circ} y^2 \ dG_{j,n}(y) \rightsquigarrow K,
\end{equation} 

\Bin under \textit{VCH[T]}.\\

\end{theorem}

\Bin Next, we have the second theorem for not necessarily centered random variables.\\

\begin{theorem} \label{main_theo2} (General Theorem for non-centered random variables) Let us adopt the notation and assumptions of Theorem \ref{main_theo1} except that the random variables are not necessarily centered. Let us denote

$$
a_{k,n}=\mathbb{E}(X_{k,n}), \ 1 \leq k\leq k(n), \ \ a_n=\sum_{k=1}^{k(n)} a_{k,n},
$$ 

$$
a_{n}^{\ast}=\sum_{j=1}^{m(n)\ell(n)} a_{j,n}, \ n\geq 1
$$ 

\Bin and

$$
b_{j,n}=\sum_{h=1}^{\ell(n)} a_{(j-1)\ell(n)+h,n}, \ 1\leq j\leq m(n).
$$

\Bin Let us denote $G_{j,n}^{\ast}(\circ)=G_{j,n}(\circ + b_{j,n})$. Then, given (\ref{C0}), (\ref{C1}) and (\ref{C2}) and assuming that \textit{UAN[T]}, \textit{BHV[T]} (or \textit{CVH[T]}) hold, we have the following results.\\

\Ni (A) If \textit{BHV[T]} holds (respectively \textit{CVH[T]} holds), if $a_n\rightarrow b$ (equivalently $a_{n}^{\ast}\rightarrow b$ if $b_n^{\ast}=\sum_{h=1}^{r(n)} a_{m(n)\ell(n)+h,n} \rightarrow 0$) and 

\begin{equation}
K_n^{\ast}(\circ)= \sum_{j=1}^{m(n)} \int_{-\infty}^{\circ} y^2 \ dG_{j,n}^{\ast}(y) \rightsquigarrow_{pre} K^{\ast}, \label{KNCB}
\end{equation} 

\Bin (resp.

\begin{equation}
K_n^{\ast}(\circ)= \sum_{j=1}^{m(n)} \int_{-\infty}^{\circ} y^2 \ dG_{j,n}^{\ast}(y) \rightsquigarrow K^{\ast}, \label{KNCC}
\end{equation} 

\Bin), then $S_n[X]\rightsquigarrow b+G$, where $G$ is \textit{idecomp} and $\psi_G(\circ)=\exp(\psi[K^\ast](\circ)$, with

\begin{equation}
\forall u \in \mathbb{R}, \ \psi[K^{\ast}](u)=\int \frac{e^{iux}-1-iux}{x^2} \ dK^{\ast}(x). \label{KC}
\end{equation}

\Bin (B) If $S_n[X] \rightsquigarrow G^{\ast}$, where $G^{\ast}$ is the \textit{cdf} of an a.s. finite random variable, then the sequence $(a_n)_{n\geq 1}$ converges to a finite number $b$ and $G^{\ast}=b+G$, with $G$ \textit{idecomp} and \eqref{KNCB} holds under the \textit{BVH[T]} or \eqref{KNCC} holds under the \textit{CVH[T]} and in both case, \eqref{KC} holds.\\
\end{theorem}

\section{General results for $G$-\textit{CLT} by-row-associated arrays} \label{sec2}

\Ni From this general description, we are going to state the overall general results for $G$-\textit{CLT} of associated random variables. As we already mentioned, existing results up to our knowledge only cover the Gaussian \textit{CLT}. \\

\Ni In the two above theorems, we need to check the conditions of validity of \textit{CLT} for by-row independent arrays on the grouped data $\{T_{j,n}, \ 1\leq j\leq m(n)\}$. However, we may get sufficient conditions on the non-grouped data $\{X_{k,n}, \ 1\leq k\leq k(n)\}$, implying them. Here are some examples. \\

\begin{corollary} \label{conv_grp_data}

\Ni Let us suppose that we have all the notation of theorems \ref{main_theo1} and \ref{main_theo2}. \\

\Ni \textbf{(1)} Suppose that the conditions (\ref{C0}), (\ref{C1}) and (\ref{C2}) hold. If \\

\Ni (i) \textit{BVH}[X] (resp. \textit{CVH}[X]) holds \\

\Ni and \\

\Ni (ii) $\ell(n) \ U(n,\epsilon/\ell(n), X)\rightarrow 0$, \\

\Ni then \textit{BVH}[T] (resp. \textit{CVH}[T]) and \textit{UAN}[T] hold and the conclusions of (1) and (2) of theorem \ref{main_theo1} are still validated. \\

\Ni \textbf{(2)}  If the Gaussian-Lynderberg condition related to $X$ is such that: 

$$
\forall \epsilon>0, \ \ \ell(n)^2 \ L_n[X](\epsilon/\ell(n))\rightarrow 0,
$$  

\Bin then the Gaussian-Lynderberg condition related to $T$ holds and we have

$$
S_{m(n)}[T]\rightsquigarrow \mathcal{N}(0,1) \ \ and \ \ B_{m(n)}[T]\rightarrow 0.
$$

\Bin Similarly, if $B_n[X]$ satisfies 

$$
\ell(n)^2 \ B_n[X]\rightarrow 0,
$$

\Bin then 

$$
S_{m(n)}[T]\rightsquigarrow \mathcal{N}(0,1) 
$$

\Bin implies that for any $\epsilon>0$,

$$
L_{m(n)}[T](\epsilon)\rightarrow 0.
$$

\end{corollary} 

\Bin \textbf{General comment}. Using the conditions on grouped and non-grouped data is a matter of situations. \\

\Ni Dealing directly with the conditions on $T_{j,n}$'s is far better and more precise. If not, we may try to establish the conditions on the non-grouped data. However, we should be aware that they are only sufficient conditions. \\

\Ni If they failing to have them, does not mean that the conditions on the grouped data fail. \\

\Ni \textbf{Elements of Proof}. These both theorems are entirely proved by the whole description in section (\ref{description}). For the proof of the Corollary \ref{conv_grp_data}, we only need to have the comparison between the condition on grouped and non-grouped data as below. \\

\Ni \textit{Conditions on the regrouped and non-regrouping by-row data}.\\

\Ni To shorten the notation, we use the following notation below

$$
X^{\prime}_{j,n,h}=X_{(j-1)\ell(n)+h,n}, \ 1\leq j\leq m(n), \ \ 1\leq h \leq \ell(n).
$$

\Bin \textbf{(1) The (\textbf{CVH}) condition}. Let us begin to see that, since $\mathbb{V}ar(S_n[X])\geq 0$ \ for any $n\geq 1$, the sequence of variances $(\mathbb{V}ar(S_n[X]))_{n\geq 1}$ converges to a some $\sigma^2>0$, otherwise, it diverges to infinity. By equation (\ref{comp_var}), we have

$$
1 = \frac{\mathbb{V}ar(S_{m(n)}[T])}{\mathbb{V}ar(S_n[X])} + \frac{\mathbb{V}ar(Y_n^{\ast})}{\mathbb{V}ar(S_n[X])} + 2\frac{\mathbb{C}ov(S_{m(n)}[T],Y_n^{\ast})}{\mathbb{V}ar(S_n[X])}
$$

\Bin and so, since $\mathbb{V}ar(Y_n^{\ast})\rightarrow 0$ (condition (\ref{C1})), and next 

$$
\frac{\mathbb{C}ov(S_{m(n)}[T],Y_n^{\ast})}{\mathbb{V}ar(S_n[X])}\rightarrow 0
$$ 

\Bin by Cauchy Schwartz inequality, we get

$$
\frac{\mathbb{V}ar(S_{m(n)}[T]}{\mathbb{V}ar(S_n[X])}\rightarrow 1.
$$

\Bin Hence

\begin{equation} \label{cond_cvh}
MV(n,X) = (1+o(1)) \ MV(m(n),T),
\end{equation}

\Bin with 

$$
MV(m(n),T) = \mathbb{V}ar(S_{m(n)}[T]). 
$$

\Bin \textbf{(2) The \textbf{UAN} condition}. We have, for any $\varepsilon>0$, the (\textbf{UAN}) condition for the sequence $\{T_{j,n}, \ 1\leq j\leq m(n)\}$ is

\begin{eqnarray*}
U(m(n),\varepsilon,T) &=& \sup_{1\leq j \leq m(n)} \mathbb{P}(|T_{j,n}|\geq \varepsilon) \\
&=& \sup_{1\leq j \leq m(n)} \mathbb{P}\left(\left|\sum_{h=1}^{\ell(n)} X^{\prime}_{j,n,h}\right|\geq \varepsilon \right) \\
&\leq& \sup_{1\leq j \leq m(n)} \mathbb{P}\left(\sum_{h=1}^{\ell(n)}\left(|X^{\prime}_{j,n,h}|>\varepsilon/\ell(n)\right)\right) \\
&\leq& \sup_{1\leq j \leq m(n)} \sum_{h=1}^{\ell(n)} \mathbb{P} \left(|X^{\prime}_{j,n,h}|>\varepsilon/\ell(n)\right) \\
&\leq& \ell(n) \sup_{1\leq j \leq m(n)}\sup_{1\leq h \leq \ell(n)} \mathbb{P} \left(|X^{\prime}_{j,n,h}|>\varepsilon/\ell(n)\right) \\
&\leq& \ell(n) \sup_{1\leq k \leq k(n)} \mathbb{P} \left(|X_{k,n}|>\varepsilon/\ell(n)\right).
\end{eqnarray*}

\Bin Hence we get for any $\epsilon>0$,

\begin{equation} \label{cond-uan}
U(m(n),\varepsilon,T)\leq \ell(n) \ U(n,\varepsilon/\ell(n),X),
\end{equation}

\Bin where we recall that 

$$
U(n,\varepsilon,X) = \sup_{1\leq k \leq k(n)} \mathbb{P} \left(|X_{k,n}|>\varepsilon \right). 
$$

\Bin \textbf{(3) A sufficient condition to obtain (\textbf{UAN}) condition}. Since for any $\varepsilon>0$,

\begin{equation*}
U(m(n),\varepsilon,T) \leq\varepsilon^{-2} \sup_{1\leq j \leq m(n)} \mathbb{E}(T_{j,n}^2), 
\end{equation*}

\Bin then the \textit{UAN} condition for the sequence $\{T_{j,n}, \ 1\leq j\leq m(n)\}$ is controlled as:

$$
B_{m(n)}[T] := \sup_{1\leq j \leq m(n)} \mathbb{E}(T_{j,n}^2) \rightarrow 0
$$

\Bin and we have

\begin{eqnarray*}
B_{m(n)}[T] &=& \sup_{1\leq j \leq m(n)} \mathbb{E} \left(\sum_{h=1}^{\ell(n)} X^{\prime}_{j,n,h}\right)^2 \\
&=& \sup_{1\leq j \leq m(n)}\mathbb{E} \left| \ell(n) \left( \sum_{h=1}^{\ell(n)} \{ X^{\prime}_{j,n,h}/\ell(n)\}\right)\right|^2\\
&=& \ell(n)^2 \sup_{1\leq j \leq m(n)} \mathbb{E} \left(\sum_{h=1}^{\ell(n)} \{ X^{\prime}_{j,n,h}/\ell(n)\}\right)^2\\
&\leq& \ell(n) \sup_{1\leq j \leq m(n)} \sum_{h=1}^{\ell(n)} \mathbb{E} \left(X^{\prime 2}_{n,j,h}\right) \ \ \ (L3) \\
&\leq& \ell(n)^2 \sup_{1\leq j \leq m(n)} \sup_{1\leq h \leq \ell(n)} \mathbb{E} \left(X^{\prime 2}_{n,j,h}\right) \\
&\leq& \ell(n)^2 \sup_{1\leq k \leq k(n)} \mathbb{E} \left(X_{k,n}^2\right),
\end{eqnarray*}

\Bin where we use in line (L3), the convexity of $\mathbb{R}_+ \ni x\mapsto x^{2}$. By denoting, for $\varepsilon>0$,

$$
B_n[X] = \sup_{1\leq k \leq k(n)} \mathbb{E} \left(X_{k,n}^2\right), 
$$

\Bin we have,

\begin{equation} \label{comp_uan}
B_{m(n)}[T]\leq \ell(n)^2 \ B_n[X].
\end{equation}

\Bin \textbf{(4) - Lyapounov condition}. For $\delta>0$, the Lyapounov condition for the sequence $\{T_{j,n}, \ 1\leq j\leq m(n)\}$ is

\begin{eqnarray*}
A_{m(n)}[T](\delta) &=& \sum_{j=1}^{m(n)} \mathbb{E} \left|T_{j,n}\right|^{2+\delta} \\
&=& \sum_{j=1}^{m(n)} \mathbb{E} \left|\sum_{h=1}^{\ell(n)} X^{\prime}_{j,n,h}\right|^{2+\delta}.
\end{eqnarray*}

\Bin Now, by using the convexity of $\mathbb{R}_+ \ni x\mapsto x^{2+\delta}$, we have
 
\begin{eqnarray*}
A_{m(n)}[T](\delta)&=& \sum_{j=1}^{m(n)} \mathbb{E} \left| \ell(n) \left( \sum_{h=1}^{\ell(n)} \{ X^{\prime}_{j,n,h}/\ell(n)\}\right)\right|^{2+\delta}\\
&\leq& \ell(n)^{1+\delta} \sum_{j=1}^{m(n)} \sum_{h=1}^{\ell(n)} \mathbb{E} \left|X^{\prime}_{j,n,h}\right|^{2+\delta}\\
&\leq& \ell(n)^{1+\delta} \sum_{k=1}^{k(n)} \left|X_{k,n}\right|^{2+\delta}.
\end{eqnarray*}

\Bin Hence, we get for any $\delta>0$,

\begin{equation}
A_{m(n)}[T](\delta) \leq \ell(n)^{1+\delta} \ A_n[X](\delta), \label{comp-lyap}
\end{equation}

\Bin where we recall that 

\begin{equation*}
A_{n}[X](\delta) = \sum_{k=1}^{k(n)} \mathbb{E} \left|X_{k,n}\right|^{2+\delta}. 
\end{equation*}

\Bin \textbf{(5) - Lynderberg-Gaussian Condition}. For $\varepsilon>0$, the Lynderberg-Gaussian Condition for the sequence $\{T_{j,n}, \ 1\leq j\leq m(n)\}$ is

\begin{equation*}
L_{m(n)}[T](\varepsilon)=\sum_{j=1}^{m(n)} \int_{\left(\left|T_{j,n}\right|\geq \varepsilon \right)} \left|T_{j,n}\right|^2 \ d\mathbb{P}.
\end{equation*}

\Bin Let us set, for $n\geq 1$, $1\leq j \leq m(n)$, $1\leq h \leq \ell(n)$,

$$
A^{\prime}_{j,n,h}=\left(\left|X^{\prime}_{j,n,h}\right|\geq \varepsilon/\ell(n)\right)
$$

\Bin and

$$
A^{\prime}_{j,n}= \bigcup_{h=1}^{\ell(n)} A^{\prime}_{j,n,h}.
$$

\Bin We have for any $\epsilon>0$,

\begin{equation*}
L_{m(n)}[T](\varepsilon)\leq \sum_{j=1}^{m(n)} \int_{A^{\prime}_{j,n}} \left|T_{j,n}\right|^2 \ d\mathbb{P}. 
\end{equation*}

\Bin Let $M_{j,n}=\max_{1\leq r\leq \ell(n)} |X^{\prime}_{j,n,r}|$. Since we have

$$
\bigcup_{r=1}^{\ell(n)} \biggr(M_{j,n}=|X^{\prime}_{j,n,r}|\biggr)=\Omega,
$$

\Bin we get

\begin{eqnarray*}
L_{m(n)}[T](\varepsilon)&\leq&  \sum_{j=1}^{m(n)} \sum_{r=1}^{\ell(n)} \int_{A^{\prime}_{j,n} \cap (M_{j,n}=|X^{\prime}_{j,n,r}|)} \left|T_{j,n}\right|^2 \ d\mathbb{P}\\
&\leq& \sum_{j=1}^{m(n)} \sum_{r=1}^{\ell(n)} \int_{A^{\prime}_{j,n} \cap (M_{j,n}=|X^{\prime}_{j,n,r}|)} \left(\ell(n) |X^{\prime}_{j,n,r}|\right)^2 \ d\mathbb{P}\\
&\leq& \ell(n)^2 \sum_{j=1}^{m(n)} \sum_{r=1}^{\ell(n)} \int_{(|X^{\prime}_{j,n,r}|\geq \varepsilon /\ell(n))} X^{\prime 2}_{j,n,r} \ d\mathbb{P}\\
&\leq& \ell(n)^2 \sum_{k=1}^{k(n)}  \int_{(|X_{k,n}|\geq \varepsilon/\ell(n))} X_{k,n}^2 \ d\mathbb{P}.
\end{eqnarray*}

\Bin Let us put for $\varepsilon>0$,

$$
L_{n}[X](\varepsilon) = \sum_{k=1}^{k(n)}  \int_{(|X_{k,n}|\geq \varepsilon)} X_{k,n}^2 \ d\mathbb{P}.
$$

\Bin Hence, we get for any $\varepsilon>0$,

\begin{equation} \label{comp_lynd_pois}
L_{m(n)}[T](\varepsilon)\leq \ell(n)^2 \ L_{n}[X](\varepsilon/\ell(n)). 
\end{equation}

\Bin We are going to specialize these results for stationary associated sequences. \\

\section{$G$-\textit{CLT} for stationary associated sequences} \label{sec3}

\Bin We will rediscover here the initial Gaussian \textit{CLT} of Newman without any other condition that $\mathbb{V}ar(S_n[X])$ converges to $\sigma^2>0$ for an infinite stationary sequence. But this result will be hardly extended to arrays. This is interesting since non-Gaussian \textit{CLT}'s usually arise for arrays as we will see it.\\

\Bin After the statement of the results, we will treat important examples of samples of Bernoulli and corrected geometric laws.\\

\Ni Let us begin by evaluating  Condition (\ref{C1}).\\

\subsection{(\ref{C1}) Automatically holds for stationary and associated data} \label{cond_C1}$ $\\

\Ni Let us suppose that each row $X_{k,n}$, \ $1\leq k\leq k(n)$ in the array $X$, is centered, square integrable, stationary and associated. \\

\Ni So, because of the stationarity, we have

$$
\mathbb{V}ar(S_n[X]) = k(n)\sigma_{1,n}^2 + 2 \sum_{j=2}^{k(n)} (k(n)-j+1) \mathbb{C}ov(X_{1,n},X_{j,n})
$$

\Bin and next suppose that 

$$
k(n)\sigma_{1,n}^2\rightarrow \lambda_1>0
$$

\Ni and

$$
2 \sum_{j=2}^{k(n)} (k(n)-j+1) \mathbb{C}ov(X_{1,n},X_{j,n})\rightarrow \lambda_2 \geq 0.
$$

\Bin We have

\begin{eqnarray*}
\mathbb{V}ar\left(Y_n^{\ast}\right) &=& \mathbb{V}ar\left(\sum_{h=1}^{r(n)} X_{m(n)\ell(n)+h,n}\right) \\
&=& \sum_{h=1}^{r(n)} \mathbb{V}ar\left(X_{m(n)\ell(n)+h,n}\right) \\
&+& \sum_{1\leq i\neq j\leq r(n)} \mathbb{C}ov\left(X_{m(n)\ell(n)+i,n}; X_{m(n)\ell(n)+j,n}\right) \\
&=& r(n) \sigma_{1,n}^2 + 2\sum_{j=2}^{r(n)} (r(n)-j+1) \ \mathbb{C}ov(X_{1,n},X_{j,n}).
\end{eqnarray*}

\Bin But

$$
r(n) \sigma_{1,n}^2 = \frac{r(n)}{k(n)} \times (k(n) \sigma_{1,n}^2)\rightarrow 0\times \lambda_1 = 0
$$

\Bin and

\begin{eqnarray*}
2\sum_{j=2}^{r(n)} (r(n)-j+1) \mathbb{C}ov(X_{1,n},X_{j,n}) &=& \frac{2 r(n)}{k(n)} \sum_{j=2}^{r(n)} \left(k(n)-\frac{k(n)}{r(n)}(j-1)\right) \mathbb{C}ov(X_{1,n},X_{j,n}) \\
&\leq& \frac{2 r(n)}{k(n)} \sum_{j=2}^{k(n)} \left(k(n)-j+1\right) \mathbb{C}ov(X_{1,n},X_{j,n}) \\
&\rightarrow& 0 \times \lambda_2 = 0.
\end{eqnarray*}

\Ni So 

$$
\mathbb{V}ar\left(Y_n^{\ast}\right)\rightarrow 0
$$
 
\Bin and next (\ref{C1}) holds.

\Bi

\subsection{Gaussian \textit{CLT} for an infinite associated and stationary sequence} $ $\\

\Ni Let us take for each $n\geq 1$,

$$
X_{k,n} = \frac{X_k}{\sqrt{n}}, \ 1\leq k\leq n=:k(n),
$$

\Ni where the $X_{k}$'s are centered, square integrable, stationary, associated and defined on the same probability space $(\Omega, \mathcal{A}, \mathbb{P})$. \\

\Ni We already know that Condition (\ref{C1}) holds (see Subsection \ref{cond_C1}, page \pageref{cond_C1}). \\

\Ni We let $\ell(n)\equiv \ell$ fix at the beginning. So

\begin{eqnarray*}
R_{m(n),\ell(n)}(u) &=& \frac{u^2}{2}\sum_{1\leq r\neq s\leq m(n)} \mathbb{C}ov(Y_{r,n},Y_{s,n}) \\
&=& \frac{u^2}{2n} \sum_{1\leq i\neq j\leq m(n)\ell(n)} \mathbb{C}ov(X_i,X_j) \\
&-& \frac{u^2}{2n} \sum_{h=1}^{m(n)} \sum_{(r,s)\in I_h^2:r\neq s} \mathbb{C}ov(X_s,X_r),
\end{eqnarray*}

\Ni with for $h\in \{1, \cdots, m(n)\}$,

$$
I_h = \{(h-1)\ell(n)+k, \ 1\leq k\leq \ell(n)\}.
$$

\Ni Hence by using the stationarity, we have

\begin{eqnarray}
R_{m(n),\ell(n)}(u) &=& \frac{u^2}{2n} \biggr(m(n)\ell(n) \mathbb{V}ar(X_1) + 2\sum_{j=2}^{m(n)\ell(n)} \left(m(n)\ell(n)-j+1\right) \mathbb{C}ov(X_1,X_j)\biggr) \notag\\
&-& \frac{u^2}{2n}\sum_{h=1}^{m(n)} \biggr(\ell(n) \mathbb{V}ar(X_1) + 2\sum_{j=2}^{\ell(n)} (\ell(n)-j+1) \mathbb{C}ov(X_1,X_j)\biggr) \notag\\
&=& \frac{u^2}{2} \frac{m(n)\ell(n)}{n}\biggr(\frac{2}{m(n)\ell(n)} \sum_{j=2}^{m(n)\ell(n)} \left(m(n)\ell(n)-j+1\right) \mathbb{C}ov(X_1,X_j) \label{expr_R} \\
&-& \frac{2}{\ell(n)} \sum_{j=2}^{\ell(n)} (\ell(n)-j+1) \mathbb{C}ov(X_1,X_j)\biggr). \notag
\end{eqnarray} 

\Bin We know that 

$$
\mathbb{V}ar(S_n[X]) = \mathbb{V}ar(X_1) + \frac{2}{n} \sum_{j=2}^{n} (n-j+1) \mathbb{C}ov(X_1,X_j),
$$

\Bin see (\cite{GH-gslo}). \\

\Ni Let us suppose that 

$$
\mathbb{V}ar(S_n[X])\rightarrow \sigma^2 = \mathbb{V}ar(X_1) + 2\sum_{j\geq 2} \mathbb{C}ov(X_1, X_j)=: \mathbb{V}ar(X_1) + \sigma_2^2<+\infty.
$$ 

\Bin [Recall that, by associativity, $\mathbb{C}ov(X_1, X_j)\geq 0$ for any $j\geq 2$]. \\

\Ni Then, from (\ref{expr_R}), we have

\begin{eqnarray*}
R_{m(n),\ell(n)}(u) &\rightarrow& \frac{u^2}{2} \left(\sigma_2^2 - \frac{2}{\ell} \sum_{j=2}^{\ell}(\ell-j+1) \ \mathbb{C}ov(X_1,X_j)\right) \\
&=& \frac{u^2}{2} B(\ell).
\end{eqnarray*}

\Bin But, by stationarity, we have

$$
S_{m(n)}[T] =_d \sum_{j=1}^{m(n)} Y_{j,n} =_d \sqrt{\frac{m(n)\ell(n)}{n}} \times \frac{1}{\sqrt{m(n)}} \sum_{j=1}^{m(n)} Z_{j,\ell},
$$

\Bin where the $Z_{j,\ell}$'s, \ $j=1,\cdots,m(n)$ \ are \textit{iid} and having the same law as

$$
Z := \frac{X_1 + \cdots + X_{\ell}}{\sqrt{\ell}}. 
$$

\Bin So 

$$
S_{m(n)}[T] =_d (1+o(1)) \ S_{m(n)}^{\ast}[T],
$$

\Bin with

$$
S_{m(n)}^{\ast}[T] = \frac{1}{\sqrt{m(n)}} \sum_{j=1}^{m(n)} Z_{j,\ell},
$$

\Bin where $Z_{1,\ell}, \cdots, Z_{m(n),\ell}$ are \textit{iid}$\sim Z$, with

$$
\mathbb{E}(Z) = 0 \ \ and \ \ \mathbb{V}ar(Z) = \mathbb{V}ar(S_{\ell}) =: \sigma_{\ell}^{2}. 
$$  

\Bin By the standard \textit{CLT} for \textit{iid} sequence,

$$
S_{m(n)}^{\ast}[T]\rightsquigarrow \mathcal{N}(0, \sigma_{\ell}^{2})
$$

\Bin and this implies that for any $u\in \mathbb{R}$, 

$$
\psi_{S_{m(n)}^{\ast}[T]}(u)\rightarrow e^{-\sigma_{\ell}^{2} u^2/2} \ as \ n\rightarrow +\infty. 
$$

\Bin Hence for any $\ell$ fixed, for any $u\in \mathbb{R}$,

\begin{eqnarray*}
\left|\psi_{S_{m(n)}[Y]}(u) - e^{-\sigma^{2}u^2/2}\right| &\leq& \left|e^{-\sigma^{2}u^2/2} - e^{-\sigma_{\ell}^{2}u^2/2}\right| \\
&+& \left|\psi_{S_{m(n)}[Y]}(u) - \psi_{S_{m(n)}^{\ast}[T]}(u)\right| + \left|\psi_{S_{m(n)}^{\ast}[T]}(u) - e^{-u^2 \sigma_{\ell}^{2}/2}\right|.
\end{eqnarray*}

\Bin So

\begin{eqnarray*}
\limsup_{n\rightarrow +\infty}\left|\psi_{S_{m(n)}[Y]}(u) - e^{-\sigma^{2}u^2/2}\right| 
&\leq& \left|e^{-\sigma_{\ell}^{2}u^2/2} - e^{-\sigma^{2}u^2/2}\right| + \frac{u^2}{2} B(\ell),
\end{eqnarray*}

\Bin for any $\ell$. Hence when we let $\ell \rightarrow +\infty$, we arrive at

$$
\limsup_{n\rightarrow +\infty}\left|\psi_{S_{m(n)}[Y]}(u) - e^{-\sigma^{2}u^2/2}\right| = 0
$$

\Bin since $\sigma_{\ell}^{2}\rightarrow \sigma^2$ and so $B(\ell)\rightarrow 0$. $\square$ \\

\Bi

\subsection{\textit{CLT} for an array with associated and stationary sequence of random variables by rows} $ $\\

\Ni Here, we consider an array
  
$$
X\equiv \biggr\{ \{X_{k,n}, \ 1 \leq k \leq k(n)\}, \ n\geq 1\biggr\}
$$

\Bin of random variables defined on the same probability space $(\Omega, \mathcal{A}, \mathbb{P})$ such that each row $\{X_{k,n}, \ 1 \leq k \leq k(n)\}$ is centered, square integrable, stationary and associated and so

$$
\mathbb{V}ar(S_n[X]) = k(n)\sigma_{1,n}^2 + 2 \sum_{j=2}^{k(n)} (k(n)-j+1) \mathbb{C}ov(X_{1,n},X_{j,n}).
$$ 

\Bin We require that 

$$
\mathbb{V}ar(S_n[X])\rightarrow \lambda>0.
$$

\Bin But this can be achieved for example if:

\begin{equation} \label{cv_1}
k(n)\sigma_{1,n}^2 \rightarrow \lambda_1
\end{equation}

\Bin and

\begin{equation} \label{cv_2}
2 \sum_{j=2}^{k(n)} (k(n)-j+1) \mathbb{C}ov(X_{1,n},X_{j,n}) \rightarrow \lambda_2,
\end{equation}

\Bin with $\lambda_1>0$ \ and \ $\lambda_2\geq 0$. \\

\Ni We already know that, with conditions (\ref{cv_1}) and (\ref{cv_2}), (\ref{C1}) holds (see Subsection \ref{cond_C1}, page \pageref{cond_C1}), but we are not able to directly show that (\ref{C2}) holds as we proved it in the last subsection, because each $S_{n}[X]$ uses here a different sequence of random variables, for any $n\geq 1$, and we do not know how the lines are related, contrary to the situation of an infinite sequence of stationary random variables as in Subsection \ref{cond_C1}. \\

\Ni We need to check that $R_{m(n),\ell(n)}(u)$ goes to zero for $u$ fixed, with 

\begin{eqnarray*} 
R_{m(n),\ell(n)}(u) &=& \frac{u^2}{2} \biggr\{2 \sum_{j=2}^{m(n)\ell(n)} \left(m(n)\ell(n)-j+1\right) \mathbb{C}ov(X_{1,n},X_{j,n}) \\ 
&-& 2 m(n) \sum_{j=2}^{\ell(n)} (\ell(n)-j+1) \mathbb{C}ov(X_{1,n},X_{j,n})\biggr\} \\
&=:& \frac{u^2}{2} \left(A_n - B_n\right).
\end{eqnarray*}  

\Bin The general conditions we may require are:

\begin{equation*}
A_n\rightarrow \lambda_2 \ and \ B_n\rightarrow \lambda_2.  
\end{equation*}

\Bin But, we have

$$
S_n[X] = S_{m(n)\ell(n)} + Y_n^{\ast}
$$

\Bin and so

$$
\mathbb{V}ar(S_n[X]) = \mathbb{V}ar(S_{m(n)\ell(n)}) + \mathbb{V}ar(Y_n^{\ast}) + 2\mathbb{C}ov(S_{m(n)\ell(n)},Y_n^{\ast}),
$$

\Bin which implies that

\begin{equation} \label{comp_var}
1 = \frac{\mathbb{V}ar(S_{m(n)\ell(n)})}{\mathbb{V}ar(S_n[X])} + \frac{\mathbb{V}ar(Y_n^{\ast})}{\mathbb{V}ar(S_n[X])} + 2\frac{\mathbb{C}ov(S_{m(n)\ell(n)},Y_n^{\ast})}{\mathbb{V}ar(S_n[X])}.
\end{equation}

\Bin Moreover, by Cauchy Schwartz inequality, we have

\begin{eqnarray*}
\mathbb{C}ov(S_{m(n)\ell(n)},Y_n^{\ast}) &\leq& \mathbb{V}ar(S_{m(n)\ell(n)})^{1/2} \ \mathbb{V}ar(Y_n^{\ast})^{1/2} \\
&\leq& \mathbb{V}ar(S_n[X])^{1/2} \ \mathbb{V}ar(Y_n^{\ast})^{1/2},
\end{eqnarray*}

\Bin where we used the association to bound $\mathbb{V}ar(S_{m(n)\ell(n)})$ by $\mathbb{V}ar(S_n[X])$. \\

\Ni Hence, by using Condition (\ref{C1}) and the hypotheses (\ref{cv_1}) and (\ref{cv_2}), we get

$$
\frac{\mathbb{V}ar(S_{m(n)\ell(n)})}{\mathbb{V}ar(S_n[X])}\rightarrow 1.
$$

\Bin Now, we have

$$
\mathbb{V}ar(S_{m(n)\ell(n)}) = m(n)\ell(n) \sigma_{1,n}^2 + A_n
$$

\Bin and since 

$$
\frac{m(n)\ell(n)}{k(n)}\rightarrow 1,
$$

\Ni we get 

$$
A_n\rightarrow \lambda_2.
$$

\Bin Hence, the only condition we require is

\begin{equation*} 
\text{(GC)} \ \ \ \ \ \ \ 2 m(n) \sum_{j=2}^{\ell(n)} (\ell(n)-j+1) \mathbb{C}ov(X_{1,n},X_{j,n}) \rightarrow \lambda_2.
\end{equation*}

\Bin Let us give two simple conditions under which (GC) holds. \\

\Ni \textbf{(PC1)} \ For $\lambda_2 = 0$,

$$
m(n)\ell(n)^2 \sup_{2\leq j\leq k(n)} \mathbb{C}ov(X_{1,n}, X_{j,n})\rightarrow 0.
$$

\Bin \textbf{(PC2)} \ For $\lambda_2>0$,

$$
\mathbb{C}ov(X_{1,n}, X_{j,n}) = \frac{\lambda_2 + o(1)}{m(n)\ell(n)^2}, \ uniformly \ in \ j\in \{2,\cdots,\ell(n)\}.
$$

\Bin Let us summarize this into the following theorem. \\

\begin{theorem}

\Ni Let

$$
X\equiv \biggr\{ \{X_{k,n}, \ 1 \leq k \leq k(n)\}, \ n\geq 1\biggr\},
$$

\Bin be an array of random variables defined on the same probability space $(\Omega, \mathcal{A}, \mathbb{P})$ such that each row $\{X_{k,n}, \ 1 \leq k \leq k(n)\}$ is centered, square integrable, stationary and associated and that

\begin{equation*} 
k(n)\sigma_{1,n}^2 \rightarrow \lambda_1
\end{equation*}

\Bin and

\begin{equation*} 
2 \sum_{j=2}^{k(n)} (k(n)-j+1) \mathbb{C}ov(X_{1,n},X_{j,n}) \rightarrow \lambda_2,
\end{equation*}

\Bin with $\lambda_1>0$ \ and \ $\lambda_2\geq 0$. \\

\Ni If one of the conditions \textbf{(PC1)} or \textbf{(PC2)} holds, then $S_n[X]$ and $S_{m(n)}[T]$ have the same weak limit law or do not.
\end{theorem}

\section{Two iconic examples} \label{sec4}

\Ni Let us extend the two classical examples of Poisson-\textit{CLT} in the independent case to associated case. \\

\begin{proposition} \label{piidBinomial}

Let $X_n\sim \mathcal{B}(n, p_n)$, $n\geq 1$, with

$$
p_n\rightarrow 0 \ \ and \ \ n p_n\rightarrow \lambda>0.
$$ 

\Bin Then 

$$
X_n\rightsquigarrow \mathcal{P}(\lambda), \ as \ n\rightarrow +\infty. \\
$$
\end{proposition}

\begin{proposition} \label{niidBinomial}

Let $X_k\sim \mathcal{N}\mathcal{B}(k, p_k)$, $k\geq 1$, with

$$
p_k\rightarrow 1 \ \ and \ \ k(1 - p_k)\rightarrow \lambda>0.
$$ 
 
\Bin Then 

$$
X_k - k\rightsquigarrow \mathcal{P}(\lambda), \ as \ k\rightarrow +\infty.
$$
\end{proposition}

\Bin Proofs of Proposition \ref{piidBinomial} and \ref{niidBinomial} can be found in \cite{GH-ips-wcia-ang}, chapter 1. \\

\Ni Now, we are going to extend these results for arrays of associated data by row. \\

\Ni Here, we suppose that for each $n\geq 1$, the elements of the row $\{X_{k,n}, \ 1\leq k\leq k(n)\}$ are associated and for any $k\in \{1, \ldots, k(n)\}$,

$$
X_{k,n}\sim \mathcal{B}(p_{k,n}), \ \ or \ \ X_{k,n}\sim \mathcal{G}(p_{k,n})-1, \ \ 0<p_{k,n}<1.
$$ 

\Bin First, we focus on the stationary case. \\

\Ni \textbf{(A)} \textbf{Stationary case with $p_{k,n}=p_n$, $1\leq k\leq k(n)$.} \\

\Ni Let us proceed for each case. \\

\Ni \textbf{(A1)} \textbf{Case of sums of associated Bernouilli laws}. We have the following result. \\

\begin{theorem} \label{conv_Bern_Pois}
Let $X$ be an array of by-row stationary associated random variables such that: \\

\Ni (i) $\forall 1\leq k\leq k(n)$,\ $X_{k,n}\sim \mathcal{B}(p_{n})$; \\

\Ni (ii) $p_n\rightarrow 0$ \  and \ $k(n) p_n\rightarrow \lambda_1$, \ as $n\rightarrow +\infty$; \\

\Ni (iii) Condition \eqref{cv_2} holds (which is justified here by \textbf{(PC1)}). \\

\Ni Then 

$$
S_n[X]\rightsquigarrow \mathcal{P}(\lambda_1).
$$

\end{theorem}

\Bin That result can be extend almost word by word for an array of corrected Geometric laws. \\

\Ni \textbf{(A2)} \textbf{Case of sums of associated corrected Geometric laws}. We have the following result. \\

\begin{theorem} \label{conv_Geom_Pois}
Let $X$ be an array of by-row stationary associated random variables such that: \\

\Ni (i) each $X_{k,n}$, $1\leq k\leq k(n)$,\ follows the corrected Geometric law of parameter $p_n$, that is $X_{k,n}$ has the same law of $Z_n-1$, where $Z_n$ is a Geometric-law of parameter $p_n$; \\

\Ni (ii) $p_n\rightarrow 1$ \ and \ $k(n)(1-p_n)\rightarrow \lambda_1>0$, \ as $n\rightarrow +\infty$; \\

\Ni (iii) Condition \eqref{cv_2} holds. \\

\Ni Then 

$$
S_n[X]\rightsquigarrow \mathcal{P}(\lambda_1).    
$$

\end{theorem}

\Bin \textbf{Proof of Theorem \ref{conv_Bern_Pois}}. The condition \eqref{cv_2} implies (\ref{C2}) holds. By using the stationarity and Conditions \eqref{cv_2}  and

$$
k(n) p_n\rightarrow \lambda_1>0,
$$

\Bin we have for $n$ large enough,

$$
\mathbb{V}ar(S_n[X]) = (1+o(1)) \ k(n) p_n + \sum_{j=2}^{k(n)} (k(n)-j+1) \ \mathbb{C}ov(X_{1,n}, X_{j,n}) \rightarrow \lambda=\lambda_1 + \lambda_2>0
$$ 

\Bin and hence Condition (\ref{C1}) automatically holds (see Subsection \ref{cond_C1}, page \pageref{cond_C1}). \\

\Ni So, for any $u\in \mathbb{R}$,

$$
\left|\psi_{S_n[X]}(u)-\psi_{S_{m(n)}[T]}(u) \right|  \rightarrow 0
$$

\Bin and

$$
\left|\psi_{S_{m(n)}[T]}(u)-\prod_{j=1}^{m(n)}\psi_{T_{j,n}}(u) \right|  \rightarrow 0.
$$

\Bin Now, we may conclude in the following way: for any $u\in \mathbb{R}$, for any $n\geq 1$,

\begin{eqnarray*}
\left|\psi_{S_n[X]}(u)-\exp\left(\lambda_1(e^{iu}-1)\right) \right| 
&\leq& \left|\psi_{S_n[X]}(u)-\psi_{S_{m(n)}[T]}(u) \right| + \left|\psi_{S_{m(n)}[T]}(u)-\prod_{j=1}^{m(n)}\psi_{T_{j,n}}(u) \right| \\
&+& \left|\prod_{j=1}^{m(n)}\psi_{T_{j,n}}(u) - \prod_{j=1}^{m(n)} \prod_{h=1}^{\ell(n)}\psi_{X^{\prime}_{h,j,n}}(u) \right| \\
&+& \left|\prod_{j=1}^{m(n)} \prod_{h=1}^{\ell(n)}\psi_{X^{\prime}_{h,j,n}}(u) - \exp\left(\lambda_1(e^{iu}-1)\right) \right| \\
&=& R_{1,n}(u) + R_{2,n}(u) + R_{3,n}(u) + R_{4,n}(u).
\end{eqnarray*}

\Bin We already know that for any $u\in \mathbb{R}$,

$$
R_{1,n}(u)\rightarrow 0 \ \ and \ \ R_{2,n}(u)\rightarrow 0.
$$

\Bin Next, we are going to use the simple following fact. \\

\Ni \textbf{Fact 1.} Let $\{(x_i,y_i),\ i\in \{1,\cdots,n\},\ n\geq 1\}$ be complex numbers such that 

$$
\forall \ i\in \{1,\cdots,n\}, \ |x_i|\leq 1 \ \ and \ \ |y_i|\leq 1.
$$

\Bin Then 

$$
\left|\prod_{i=1}^{n} x_i - \prod_{i=1}^{n} y_i\right|\leq \sum_{i=1}^{n} |x_i-y_i|.
$$

\Bin Moreover, for any $u\in \mathbb{R}$,

\begin{eqnarray}
R_{3,n}(u) &=& \left|\prod_{j=1}^{m(n)}\psi_{T_{j,n}}(u) - \prod_{j=1}^{m(n)} \prod_{h=1}^{\ell(n)}\psi_{X^{\prime}_{h,j,n}}(u) \right| \notag \\
&\leq& \sum_{j=1}^{m(n)} \left|\psi_{T_{j,n}}(u) - \prod_{h=1}^{\ell(n)}\psi_{X^{\prime}_{h,j,n}}(u) \right| \ \ \ (L2) \notag \\
&\leq& \sum_{j=1}^{m(n)} \frac{u^2}{2}\sum_{1\leq r \neq s\leq \ell(n)} \mathbb{C}ov(X^{\prime}_{r,j,n}, \ X^{\prime}_{s,j,n}) \ \ \ (L3) \notag \\
&\leq& \frac{u^2}{2} m(n) \sup_{1\leq j\leq m(n)}\sum_{1\leq r \neq s\leq \ell(n)} \mathbb{C}ov(X^{\prime}_{r,j,n}, \ X^{\prime}_{s,j,n}),  \label{R3n}
\end{eqnarray}

\Bin where we use Fact 1 in the line (L2) and Newman's Lemma (See \cite{GH-newmanwright, GH-newmanwright2}) in the line (L3). \\ 

\Ni Hence, we use the stationarity in formula \eqref{R3n} to get, for any $u\in \mathbb{R}$,  

\begin{eqnarray*}
R_{3,n}(u) &\leq& u^2 m(n) \sup_{1\leq j\leq m(n)} \sum_{h=2}^{\ell(n)} (\ell(n)-h+1) \ \mathbb{C}ov(X_{1,n},\ X_{h,n}) \\
&\leq& u^2 m(n)\ell(n)^2 \sup_{2\leq h\leq \ell(n)} \ \mathbb{C}ov(X_{1,n},\ X_{h,n}),
\end{eqnarray*}

\Bin which converges to zero by \textbf{(PC1)}. \\

\Ni Finally, for any $u\in \mathbb{R}$,

\begin{eqnarray*}
C_n(u) &:=& \prod_{j=1}^{m(n)} \prod_{h=1}^{\ell(n)}\psi_{X^{\prime}_{h,j,n}}(u) \\
&=& \left((1-p_n) + p_n e^{iu}\right)^{m(n)\ell(n)} \\
&=& \left(\left((1-p_n) + p_n e^{iu}\right)^{k(n)}\right)^{m(n)\ell(n)/k(n)}.
\end{eqnarray*}

\Bin Here, by Proposition \ref{piidBinomial}, we have the asymptotic law of 

$$
Z_n \sim \mathcal{B}(k(n), p_n),
$$

\Bin with 

$$
p_n\rightarrow 0 \ \ and \ \ k(n)p_n \rightarrow \lambda_1
$$

\Bin to Poisson law $\mathcal{P}(\lambda)$. \\

\Ni Then for any $u\in \mathbb{R}$,

$$
\left((1-p_n) + p_n e^{iu}\right)^{k(n)}\rightarrow \exp\left(\lambda_1(e^{iu}-1)\right)
$$

\Bin and next

$$
C_n(u) \rightarrow \exp\left(\lambda_1(e^{iu}-1)\right)
$$

\Bin since 

$$
\frac{m(n)\ell(n)}{k(n)}\rightarrow 1.
$$

\Bin The proof is over.  $\square$ \\

\Bin \textbf{Proof of Theorem \ref{conv_Geom_Pois}}. We are going to follow the lines of the proof of Theorem \ref{conv_Bern_Pois}. Let us put $q_n:=1-p_n$. So, by using the stationarity and condition (ii), we have for $n$ large enough,

\begin{eqnarray*}
\mathbb{V}ar(S_n[X]) &=& k(n)\ \frac{q_n}{p_n^2} + \sum_{j=2}^{k(n)} (k(n)-j+1) \ \mathbb{C}ov(X_{1,n}, X_{j,n}) \\
&=& (1+o(1))\ k(n) q_n + \sum_{j=2}^{k(n)} (k(n)-j+1) \ \mathbb{C}ov(X_{1,n}, X_{j,n})
\end{eqnarray*}

\Bin and next by Conditions (ii) and (iii),

$$
\mathbb{V}ar(S_n[X])\rightarrow \lambda_1 + \lambda_2=:\lambda>0.
$$

\Bin Hence, Conditions (\ref{C1}) and (\ref{C2})  hold (see Subsection \ref{cond_C1}, page \pageref{cond_C1} and \eqref{cv_2}) and next for any $u\in \mathbb{R}$,

$$
\left|\psi_{S_n[X]}(u)-\psi_{S_{m(n)}[T]}(u) \right|  \rightarrow 0
$$

\Bin and

$$
\left|\psi_{S_{m(n)}[T]}(u)-\prod_{j=1}^{m(n)}\psi_{T_{j,n}}(u) \right|  \rightarrow 0.
$$

\Bin Now, we are going to use the same method as the proof of Theorem \ref{conv_Bern_Pois} to show that: for any $u\in \mathbb{R}$, for any $n\geq 1$,

$$
\left|\psi_{S_n[X]}(u)-\exp\left(\lambda_1(e^{iu}-1)\right) \right| \leq R_{1,n}(u) + R_{2,n}(u) + R_{3,n}(u) + R_{4,n}(u),
$$

\Bin with the $R_{q,n}(u)$, \ $q\in \{1,2,3,4\}$, are already defined in the proof of Theorem \ref{conv_Bern_Pois} and for any $u\in \mathbb{R}$,

$$
R_{q,n}(u)\rightarrow 0 \  for \ q\in \{1,2,3\}.
$$

\Bin Finally, for any $u\in \mathbb{R}$,

\begin{eqnarray*}
C_n(u) &:=& \prod_{j=1}^{m(n)} \prod_{h=1}^{\ell(n)}\psi_{X^{\prime}_{h,j,n}}(u) \\
&=& \left(\frac{p_n}{1-q_n e^{iu}}\right)^{m(n)\ell(n)} \\
&=& \left(\left(\frac{p_n}{1-q_n e^{iu}}\right)^{k(n)}\right)^{m(n)\ell(n)/k(n)}.
\end{eqnarray*}

\Bin Let us conclude as follows: let 

$$
Z_n \sim \mathcal{NB}(k(n), p_n),
$$

\Bin with here

$$
p_n\rightarrow 1 \ \ and \ \ k(n)q_n \rightarrow \lambda_1.
$$

\Bin Then, by Proposition \ref{niidBinomial},

$$
Z_{n}-n \rightsquigarrow \mathcal{P}(\lambda)
$$

\Bin and this implies that, for any $u\in \mathbb{R}$,

$$
\left(\frac{p_n}{1-q_n e^{iu}}\right)^{k(n)} = \Phi_{Z_{n}-n}(u) \rightarrow \exp\left(\lambda_1(e^{iu}-1)\right)
$$

\Bin and next

$$
C_n(u) \rightarrow \exp\left(\lambda_1(e^{iu}-1)\right)
$$

\Bin since 

$$
\frac{m(n)\ell(n)}{k(n)}\rightarrow 1.
$$

\Bin The proof is over.  $\square$ \\

\Ni \textbf{(B)} \textbf{General case.} We are going to generalize these results above for an array non necessarily stationary by row. \\

\Ni \textbf{(B1)} \textbf{Case of sums of associated Bernouilli laws}. We have the following result. \\

\begin{theorem} \label{conv_Bern_Pois_Gen}
Let $X$ be an array of by-row associated random variables such that: \\

\Ni (1) $\forall 1\leq k\leq k(n)$,

$$
X_{k,n}\sim \mathcal{B}(p_{k,n});
$$

\Bin (2) as $n\rightarrow +\infty$, 

$$
(a) \ \ \overline{p}_n := \sup_{1\leq k\leq k(n)} p_{k,n}\rightarrow 0,
$$ 

$$
(b) \ \ \sum_{k=1}^{k(n)} p_{k,n} \rightarrow \lambda>0,
$$  

$$
(c) \ \ k(n) \ \overline{p}_n^2\rightarrow 0
$$ 

\Ni and 

$$
(d) \ \ k(n)^2 \ \sup_{1\leq r\neq s\leq k(n)} \mathbb{C}ov(X_{r,n}, \ X_{s,n}) \rightarrow 0;
$$

\Bin (3)  as $n\rightarrow +\infty$, 

$$
\sum_{h=1}^{r(n)} p_{m(n)\ell(n)+h,n}\rightarrow 0;
$$

\Bin (4)  

$$
(a) \ \ \lim_{n\rightarrow +\infty} \ m(n) \sup_{1\leq j\leq m(n)} \sum_{1\leq r \neq s\leq \ell(n)} \mathbb{C}ov(X^{\prime}_{r,j,n}, \ X^{\prime}_{s,j,n}) = 0
$$ 

\Bin and   

$$
(b) \ \ \lim_{n\rightarrow +\infty}\sum_{1\leq r \neq s\leq r(n)} \mathbb{C}ov(X_{m(n)\ell(n)+r,n}, \ X_{m(n)\ell(n)+s,n})\rightarrow 0.
$$

\Bin Then 

$$
S_n[X]\rightsquigarrow \mathcal{P}(\lambda).
$$

\end{theorem}

\Bin \textbf{Proof of Theorem \ref{conv_Bern_Pois_Gen}.} Throughout this proof, the notation $\ell_{k,n}=\overline{o}_n(1)$, for $k$ ranging over some set $I_n$ means that the sequence $\ell_{k,n}$ goes to zero as $n\rightarrow +\infty$ uniformly in $k \in I_n$. So Assumption (2a) means that

$$
p_{k,n}=\overline{o}_n(1) \ \ and \ \ q_{k,n}=1-p_{k,n}=1+\overline{o}_n(1).
$$

\Bin We have the convergence of variances since by (2a),

$$
\mathbb{V}ar(S_n[X]) = (1+\overline{o}_n(1)) \ \sum_{k=1}^{k(n)} p_{k,n} + \sum_{1\leq r \neq s\leq k(n)} \mathbb{C}ov(X_{r,n}, \ X_{s,n}),
$$

\Bin which converges by (2b) and (2d). \\

\Ni Next, we can see that Condition (\ref{C1}) holds by (3) and (4b). Also, Condition \eqref{C2} holds by (2d). Indeed, that last Condition is controlled by 

$$
\sum_{1\leq r\neq s\leq m(n)} \mathbb{C}ov(T_{r,n}, \ T_{s,n})\rightarrow 0
$$ 

\Bin and that

\begin{eqnarray*}
\sum_{1\leq r\neq s\leq m(n)} \mathbb{C}ov(T_{r,n}, \ T_{s,n}) &=& \sum_{1\leq r\neq s\leq m(n)\ell(n)} \mathbb{C}ov(X_{r,n}, \ X_{s,n}) \\
&-& \sum_{j=1}^{m(n)} \sum_{1\leq r \neq s\leq \ell(n)} \mathbb{C}ov(X^{\prime}_{r,j,n}, \ X^{\prime}_{s,j,n}) \\
&\leq& \sum_{1\leq r\neq s\leq k(n)} \mathbb{C}ov(X_{r,n}, \ X_{s,n}) \ \ \ \ (L2) \\
&\leq& k(n)^2 \ \sup_{1\leq r\neq s\leq k(n)} \mathbb{C}ov(X_{r,n}, \ X_{s,n}),
\end{eqnarray*}

\Bin where the line (L2) is justified by the non-negativity of covariances for associated data. \\

\Ni Hence, we are going to use the same method as in the Proof of Theorem \ref{conv_Bern_Pois} to conclude by showing 

$$
\forall u\in \mathbb{R}, \ R_{q,n}(u)\rightarrow 0, \ as \ n\rightarrow +\infty \ \ for \ \ q\in \{1,2,3,4\}.
$$

\Bin But it is clear that for any $u\in \mathbb{R}$,

$$
R_{q,n}(u)\rightarrow 0, \ as \ n\rightarrow +\infty \ \ for \ \ q\in \{1,2\}.
$$

\Bin Moreover, we easily see by \eqref{R3n} and the hypothesis (4a) that, for any $u\in \mathbb{R}$,

\begin{equation} \label{cond_R3n_bern_non_assoc}
R_{3,n}(u)\rightarrow 0, \ as \ n\rightarrow +\infty.
\end{equation}

\Bin Hence, we should show that for any $u\in \mathbb{R}$,

$$
\overline{\phi}_n (u) := \prod_{j=1}^{m(n)} \prod_{h=1}^{\ell(n)}\psi_{X^{\prime}_{h,j,n}}(u)\rightarrow \exp\left(\lambda(e^{iu}-1)\right)
$$

\Bin to conclude. But, for $u$ fixed, we have

$$
\overline{\phi}_n (u) = \prod_{j=1}^{m(n)} \prod_{h=1}^{\ell(n)} \left(1 + p_{(j-1)\ell(n)+h,n}(e^{iu} - 1)\right)
$$

\Bin and hence 

\begin{eqnarray*}
\log\overline{\phi}_n (u) &=& \sum_{j=1}^{m(n)} \sum_{h=1}^{\ell(n)} \left(p_{(j-1)\ell(n)+h,n}(e^{iu} - 1) + O(\overline{p}_n^2)\right) \\
&=& \left(\sum_{k=1}^{m(n)\ell(n)} p_{k,n}\right)(e^{iu}-1) + O\left(k(n) \ \overline{p}_n^2\right).
\end{eqnarray*}

\Bin By (2b) and (3),

$$
\sum_{k=1}^{m(n)\ell(n)} p_{k,n} \rightarrow \lambda
$$

\Bin and hence, we conclude by (2c) that,
$$
\log\overline{\phi}_n (u)\rightarrow \lambda(e^{iu}-1),
$$

\Bin that is,

$$
\overline{\phi}_n (u)\rightarrow \exp\left(\lambda(e^{iu}-1)\right) = \phi_{\mathcal{P}(\lambda)}(u). \ \ \square
$$

\Bin \textbf{Remark 1.} The assumption (2c) is very reasonable, since for the stationary case, we have

$$
\overline{p}_n = p_n
$$ 

\Ni and since

$$
k(n) p_n\rightarrow \lambda,
$$

\Ni we have

$$
k(n)\ \overline{p}_n^2 = \frac{(k(n) p_n)^2}{k(n)} = \frac{\lambda^2(1+o(1))}{k(n)} \rightarrow 0.
$$

\Bin \textbf{(B2)} \textbf{Case of sums of associated corrected Geometric laws}. We have the following result. \\

\begin{theorem} \label{conv_Geom_Pois_Gen}
Let $X$ be an array of by-row associated random variables such that: \\

\Ni (1) $\forall 1\leq k\leq k(n)$,

$$
X_{k,n}=_{d} Z_{k,n} -1, \ \ 	Z_{k,n}\sim \mathcal{G}(p_{k,n});
$$

\Bin (2) for $q_{k,n} = 1 - p_{k,n}$, \ as $n\rightarrow +\infty$, 

$$
(a) \ \ \overline{q}_n := \sup_{1\leq k\leq k(n)} q_{k,n}\rightarrow 0,
$$ 

\Ni and

$$
(b) \ \ k(n) \ \overline{q}_n^2 \rightarrow 0;
$$ 

\Bin (3)  as $n\rightarrow +\infty$, 

$$
(a) \ \ \sum_{k=1}^{k(n)} q_{k,n}\rightarrow \lambda,
$$

$$
(b) \ \ \sum_{h=1}^{r(n)} q_{m(n)\ell(n)+h,n}\rightarrow 0
$$

\Bin and

$$
(c) \ \ k(n)^2 \ \sup_{1\leq r\neq s\leq k(n)} \mathbb{C}ov(X_{r,n}, \ X_{s,n}) \rightarrow 0;
$$

\Bin (4)  

$$
(a) \ \ \lim_{n\rightarrow +\infty} \ m(n) \sup_{1\leq j\leq m(n)} \sum_{1\leq r \neq s\leq \ell(n)} \mathbb{C}ov(X^{\prime}_{r,j,n}, \ X^{\prime}_{s,j,n}) = 0
$$ 

\Bin and   

$$
(b) \ \ \lim_{n\rightarrow +\infty}\sum_{1\leq r \neq s\leq r(n)} \mathbb{C}ov(X_{m(n)\ell(n)+r,n}, \ X_{m(n)\ell(n)+s,n})\rightarrow 0.
$$

\Bin Then 

$$
S_n[X]\rightsquigarrow \mathcal{P}(\lambda).
$$

\end{theorem}

\Bin \textbf{Remark 2.} The assumption (2b) is very reasonable, since for the stationary case, we have

$$
\overline{q}_n = q_n
$$ 

\Ni and since

$$
k(n) q_n\rightarrow \lambda,
$$

\Ni we have

$$
k(n)\ \overline{q}_n^2 = \frac{(k(n) q_n)^2}{k(n)} = \frac{\lambda^2(1+o(1))}{k(n)} \rightarrow 0.
$$

\Bin \textbf{Proof of Theorem \ref{conv_Geom_Pois_Gen}.} We are going to follow the lines of the proof of Theorem \ref{conv_Bern_Pois_Gen}. Under the hypotheses, we can easily show that Conditions (\ref{C1}) and (\ref{C2}) hold and next for any $u\in \mathbb{R}$,

$$
R_{q,n}(u)\rightarrow 0, \ as \ n\rightarrow +\infty \ \ for \ \ q\in \{1,2,3\}.
$$

\Bin Hence, we should show that for any $u\in \mathbb{R}$,

$$
\overline{\phi}_n (u) := \prod_{j=1}^{m(n)} \prod_{h=1}^{\ell(n)}\psi_{X^{\prime}_{h,j,n}}(u)\rightarrow \exp\left(\lambda(e^{iu}-1)\right)
$$

\Bin to conclude. But, for $u$ fixed, we have

$$
\overline{\phi}_n (u) = \prod_{j=1}^{m(n)} \prod_{h=1}^{\ell(n)} \left(\frac{p_{(j-1)\ell(n)+h,n}}{1-q_{(j-1)\ell(n)+h,n}e^{iu}}\right)
$$

\Bin and hence 

\begin{eqnarray*}
\log\overline{\phi}_n (u) &=& \sum_{j=1}^{m(n)} \sum_{h=1}^{\ell(n)} \left(\log(1-q_{(j-1)\ell(n)+h,n}) - \log(1-q_{(j-1)\ell(n)+h,n}e^{iu})\right)  \\ 
&=& \sum_{j=1}^{m(n)} \sum_{h=1}^{\ell(n)} \left(-q_{(j-1)\ell(n)+h,n} + q_{(j-1)\ell(n)+h,n}e^{iu} + O(\overline{q}_n^2)\right) \\
&=& \left(\sum_{k=1}^{m(n)\ell(n)} q_{k,n}\right)(e^{iu}-1) + O\left(k(n) \ \overline{q}_n^2\right).
\end{eqnarray*}

\Bin By (3a) and (3b),

$$
\sum_{k=1}^{m(n)\ell(n)} q_{k,n} \rightarrow \lambda
$$

\Bin and hence, we conclude by (2b) that,

$$
\log\overline{\phi}_n (u)\rightarrow \lambda(e^{iu}-1),
$$

\Bin that is,

$$
\overline{\phi}_n (u)\rightarrow \exp\left(\lambda(e^{iu}-1)\right) = \phi_{\mathcal{P}(\lambda)}(u). \ \ \square
$$

\section{conclusion} \label{sec5}

\Ni We have obtained a general \textit{CLT} theory of stationary and non-stationary associated arrays in the frame of Newman approximation under specific conditions.
This theory can be improved by tring to loosen the conditions. However, we think they might be optimal since we closely work on the borders of the Newman's frame. The asymptotic theory on associated data outside that frame is still a major concern and untill now, there is not a single result of this.\\

\Ni Now, it is natural to extend our results to functional limit laws (invariance principles) with a deterministic or a random horizon. This will be the main topic of a coming paper.


\begin{thebibliography}{99}

\bibitem[Adekpejou \textit{et al.} (2021)]{akim21} Adekpejou A., Traor\'e C.M.M., Lo G.S. and Niang A.B and Lo G.S. (2020).
Central limit theorems for associated possibly moving partial sums and application to the non-stationary invariance principle. arXiv:2108.10906.

\bibitem[Bulinski and Shashkin (2007)]{GH-bulinski2007} Bulinski A. and Shashkin A.(2007). Limit theorems for
associated random fields and related systems. World Scientific Publishing,
Singapore.


\bibitem[Burton \textit{et al.} (1986)]{burton} Burton, R.M., Dabrowski, A.R. and Dehling, H. (1986). An
invariance principle for weakly associated random variables, \textit{Stoch.
Proc. Appl.}, 23, 301-306.

\bibitem[Cox and Grimmett (1984)]{cox} Cox, J.T. and Grimmett, G. (1984) Central limit theorems for associated random variables and the percolation model, \textit{Ann. Probab.}, 12, 514-528.

\bibitem[Dabrowski and Dehling(1988)]{dabro} Dabrowski, A.R. and Dehling, H. (1988). A Berry-Esseen
theorem and a functional law of the iterated logarithm for weakly associated
random variables, \textit{Stochastic Process. Appl.}, 30, 247-289.

\bibitem[Esary \textit{et al.} (1967)]{esary} Esary, J., Proschan, F. and Walkup, D.(1967). Association of
random variables with application. \textit{Ann. Math Statist.}, 38.

\bibitem[Fortuin \textit{et al.} (1968)]{fortuin} Fortuin, C., Kastelyn, P. and Ginibre, J.(1971).
Correlation inequalities on some partially ordered sets. \textit{Comm. Math.
Phys.}, 22, 89-103

\bibitem[Joazhu Daley (2002)]{jiazhu} Jiazhu, P.(2002). Tail dependence of random variables from
ARCH and heavy-tailed bilinear models. \textit{Sciences in China}, 45 (6),
Ser. A, 749-760.

\bibitem[Lo (2018)]{GH-ips-wcia-ang} Lo, G.S.(2018). Weak Convergence (IA). Sequences of
random vectors. SPAS Books Series.(2016). Doi : 10.16929/sbs/2016.0001.

\bibitem[Lo \textit{et al.} and Lo (2018)]{loSanMbarka18} G. S. Lo, A. M. Fall and H. Sangaré (2018). A Central limit Theorem of dependent sums of standard exponential functionals motivated by extreme value theory. Afrika Statistika. Vol. 13 (3), 1795-1822. DOI: http://dx.doi.org/10.16929/as/1795.134 https://projecteuclid.org/euclid.as/1544583716.

\bibitem[Lo\`eve (1977)]{loeve} Lo\`{e}ve, M.(1977). \textit{Probability Theory I}.
Springer-Verlag. New-York.

\bibitem[Louhichi (2000)]{GH-louhichi} Louhichi, S.(2000). Weak convergence for empirical
processes of associated sequences. Ann. Inst. Henri Poincar\'e,
Probabilit\'es et Statistiques 36 (\textbf{5}), pp. 547--567.

\bibitem[Newman (1980)]{newman} Newman C.M. (1980) Normal fluctuations and the FKG
inequalities. \textit{Comm. Math. Phys.} 74, 119-128.

\bibitem[Newman and Wright  (1981)]{newmanwright} Newman, C.M and Wright, A.L.(1981). An invariance
principle for certain dependent sequences. \textit{Ann. probab.}, \textbf{9}
(4), 671-675.

\bibitem[Newman and Wright (1982)]{newmanwright2} Newman, C.M and Wright, A.L.(1982). Associated
random variables and martingale inequalities. Z. Wahrscheinlichkeitstheorie
verw. Gebiete 59, 361-371.

\bibitem[Newman and Wright (1981)]{GH-newmanwright} Newman, C.M and Wright, A.L.(1981). An invariance
principle for certain dependent sequences. \textit{Ann. probab.}, \textbf{9}
(4), 671-675.

\bibitem[Newman and Wright (1982)]{GH-newmanwright2} Newman, C.M and Wright, A.L.(1982). Associated
random variables and martingale inequalities. Z. Wahrscheinlichkeitstheorie
verw. Gebiete 59, 361-371.

\bibitem[Oliveira (2012)]{GH-paulo} Oliveira, P.E.(2012). \textit{Asymptotics for Associated
Random Variables.} DOI 10.1007/978-3-642-25532-8, Springer-Verlag Berlin Heidelberg.

\bibitem[Prakasa Rao (2012)]{rao} Prakasa Rao, B. L. S.(2012). \textit{Associated sequences,
Demimartingales and Nonparametric Inference.Probability and its applications}%
. Springer Basel Doredrecht, Heidelberg, London, New York.

\bibitem[Sangar\'e and Lo (2015)]{sang15} A general strong law of large numbers and applications to associated sequences and to extreme value theory, Annales Mathematicae et Informaticae, 45 (2015) pp. 111–132, http://ami.ektf.hu.

\bibitem[Sanghar\'e and Lo (2016)]{GH-gslo} Sangar\'e, H. and Lo, G. S.(2016) A Review on asymptotic normality of sums of associated random variables. \textit{Afrika Statistika}, , 11 (\textbf{1}), pp.855-867. Doi : 10.16929/as/2016.855.79. Arxiv 1405.4316.

\bibitem[Sangar\'e and Lo (2018)]{sanglo28} Sangaré H. and Lo, G. S. (2018). General Central Limit Theorems for Associated Sequences. A Collection of Papers in Mathematics and Related Sciences. Spas Editions, Euclid Series Book, pp. 289-308. Doi: http://dx.doi.org/10.16929/sbs/2018.100-04-01, \textit{url}: https://projecteuclid.org/euclid.spaseds/1569509478.

\bibitem[Sangar\'e \textit{et al.} (2020)]{sang20} Sangaré H., Lo, G. S. and Traoré M.C.M. (2020).  Arbitrary functional Glivenko-Cantelli classes and applications to different types of dependence. Far East Journal of Theoretical Statistics. Vol. 60 (1-2), 41-62. ISSN: 0972-0863. 
http://dx.doi.org/10.17654/TS060020041. 

\bibitem[Traor\'e \textit{et al.} (2016)]{traore16} Traor\'e C.M.M., Lo G.S. and Niang A.B (2016) Invariance principles for random sums of non-stationary independent and associated random variables. arXiv:1610.02700

\bibitem[Yu (1968)]{yu93} Yu, H.(1993). A Gkivenko-Cantelli lemma and weak convergence
for empirical processes of associated sequences. Probab. Theory Related
Fields 95, 357-370.
\end{thebibliography}
\end{document}